\documentclass{amsart}

\usepackage{a4,xy,amssymb,diagrams,times,array,color}
\xyoption{poly}
\xyoption{arc}
\xyoption{knot}
\xyoption{curve}
\xyoption{arrow}
\xyoption{all}

\makeindex
 
\newcommand{\N}{\mathbb{N}}

\newcommand{\vtx}[1]{*+[o][F-]{\scriptscriptstyle #1}}

\newcolumntype{C}{>{$}c<{$}}

\newcommand{\aar}[2]{ \ar@{-}|(.25){{\bf #1}}|{\bullet}|(.75){{\bf  #2}}}

\title{braid group $B_3$ irreducibles \\ - a DIY guide - }
\author{Lieven Le Bruyn} 
\address{Department of Mathematics, University of Antwerp \\ 
 Middelheimlaan 1, B-2020 Antwerp (Belgium) \\ {\tt lieven.lebruyn@ua.ac.be}}

\begin{document}
\sloppy
 
\def\ldb{\mathopen{\{\!\!\{}} \def\rdb{\mathclose{\}\!\!\}}}


 \begin{abstract} This note tells you how to construct a $k(n)$-dimensional family of (isomorphism classes of) irreducible representations of dimension $n$ for the three string braid group $B_3$, where $k(n)$ is an admissible function of your choosing; for example take $k(n) = [ \tfrac{n}{2} ] +1$ as in \cite{brol1} and \cite{brol2}.
\end{abstract}

\maketitle

\noindent
{\bf (step 1)} {\bf Learn the basics. } The three string braid group $B_3$ is the group $\langle \sigma_1,\sigma_2 | \sigma_1 \sigma_2 \sigma_1 = \sigma_2 \sigma_1 \sigma_2 \rangle$ and its center is cyclic with generator $c = (\sigma_1 \sigma_2)^3 = (\sigma_1 \sigma_2 \sigma_1)^2$. The quotient group
\[
B_3 / \langle c \rangle = \langle u,v | u^2=v^3=e \rangle \simeq C_2 \ast C_3 \simeq \Gamma_0 \]
is the modular group $PSL_2(\mathbb{Z})$ where $u$ and $v$ are the images of $\sigma_1 \sigma_2$ resp. $\sigma_1 \sigma_2 \sigma_1$.

By Schur's lemma, the central element $c$ acts as $\lambda I_n$ (where $\lambda \in \mathbb{C}^*$) on any $n$-dimensional irreducible $B_3$-representation. Hence, it is enough to construct a $k(n)-1$-dimensional family of $n$-dimensional irreducible representations of the modular group $\Gamma_0$.

If $V$ is an $n$-dimensional $\Gamma_0$ representation, we can decompose it into eigenspaces for the action of $C_2 = \langle u \rangle$ and $C_3 = \langle v \rangle$ :
\[
V_1 \oplus V_2 = V \downarrow_{C_2} = V = V \downarrow_{C_3} = W_1 \oplus W_2 \oplus W_3 \]
If the dimension of $V_i$ is $a_i$ and that of $W_j$ is $b_j$, we say that $V$ is a $\Gamma_0$-representation of {\em dimension vector} $\alpha = (a_1,a_2;b_1,b_2,b_3)$. Choosing a basis $B_1$ of $V$ wrt. the decomposition $V_1 \oplus V_2$ and a basis $B_2$ wrt. $W_1 \oplus W_2 \oplus W_3$, we can view the basechange matrix $B_1 \rTo B_2$ as an $\alpha$-dimensional representation $V_Q$ of the quiver $Q$
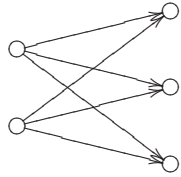
\[
\xymatrix@=.3cm{
& & & & \vtx{} \\
\vtx{} \ar[rrrru] \ar[rrrrd] \ar[rrrrddd] & & & & \\
& & & & \vtx{} \\
\vtx{} \ar[rrrru] \ar[rrrruuu] \ar[rrrrd] & & & & \\
& & & &  \vtx{} }
\]
Bruce Westbury \cite{Westbury} has shown that $V$ is an irreducible $\Gamma_0$-representation if and only if $V_Q$ is a {\em $\theta$-stable} $Q$-representation where $\theta = (-1,-1;1,1,1)$ and that the two notions of isomorphism coincide. The {\em Euler-form} $\chi_Q$ of the quiver $Q$ is the bilinear form on $\mathbb{Z}^{\oplus 5}$ determined by the matrix
\[
\begin{bmatrix} 1 & 0 & -1 & -1 & -1 \\ 0 & 1 & -1 & -1 & -1 \\ 0 & 0 & 1 & 0 & 0 \\ 0 & 0 & 0 & 1 & 0 \\ 0 & 0 & 0 & 0 & 1 \end{bmatrix} \]
Westbury also showed that if there exists a $\theta$-stable $\alpha$-dimensional $Q$-representation, then there is an $1 - \chi_Q(\alpha,\alpha)$ dimensional family of isomorphism classes of such representations (and a Zariski open subset of them will correspond to isomorphism classes of irreducible $\Gamma_0$-representations). Hence, an {\em admissible} function $k(n)$ is one such that for all $n$ we have $k(n) \leq 2 - \chi_Q(\alpha_n,\alpha_n)$ for a dimension vector $\alpha_n=(a_1,a_2;b_1,b_2,b_3)$ such that $n=a_1+a_2$ and there exists a $\theta$-stable $\alpha_n$-dimensional $Q$-representation. Note that Aidan Schofield \cite{Schofield} gave an inductive procedure to determine the dimension vectors of stable representations.

\par \vskip 4mm

\noindent
{\bf (step 2) Choose  known non-isomorphic $\Gamma_0$-irreducibles } and their corresponding $\theta$-stable $Q$-representations $\{ V_i~:~i \in I \}$. Here are some obvious choices : using the foregoing and standard quiverology, there are $6$ irreducible $1$-dimensional $\Gamma_0$-representations $S_{ij}$ and there are $3$ one-parameter families of $2$-dimensional simple $\Gamma_0$-representations $T_i(\lambda)$. Below the corresponding $Q$-representations for $S_{21}$ and $T_2(\lambda)$ (the other cases are similar)
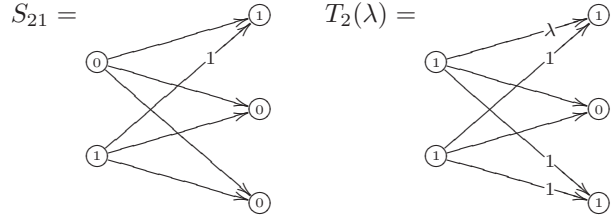
\[
S_{21} = \xymatrix@=.3cm{
& & & & \vtx{1} \\
\vtx{0} \ar[rrrru] \ar[rrrrd] \ar[rrrrddd] & & & & \\
& & & & \vtx{0} \\
\vtx{1} \ar[rrrru] \ar[rrrruuu]|(0.7){1} \ar[rrrrd] & & & & \\
& & & &  \vtx{0} } \qquad 
T_2(\lambda) = \xymatrix@=.3cm{
& & & & \vtx{1} \\
\vtx{1} \ar[rrrru]|(0.7){\lambda} \ar[rrrrd] \ar[rrrrddd]|(0.7){1} & & & & \\
& & & & \vtx{0} \\
\vtx{1} \ar[rrrru] \ar[rrrruuu]|(0.7){1} \ar[rrrrd]|(0.7){1} & & & & \\
& & & &  \vtx{1} }
\]
More interesting choices are the $Q$-representations corresponding to irreducible continuous representations of $\hat{\Gamma}_0$, the profinite completion of the modular group. For example, a simple factor of the monodromy representation associated to a dessin d'enfant or an irreducible representation of a finite group generated by an order two and an order three element, for example the monster group $\mathbb{M}$. Pick your favourite collection of non-isomorphic $\{ V_i \}$.

\par \vskip 4mm

\noindent
{\bf (step 3) Compute the local quiver} of the collection $\{ V_i~:~i \in I \}$ as in e.g. \cite{AdriLieven}. That is, we make a new quiver $\Delta$ having one vertex $v_i$ for every $V_i$. If $\alpha_i$ is the dimension vector of the $\theta$-stable $Q$-representation determined by $V_i$, then there are $1 - \chi_Q(\alpha_i,\alpha_i)$ loops in vertex $v_i$ in $\Delta$ and there are exactly $-\chi_Q(\alpha_i,\alpha_j)$ oriented arrows starting in vertex $v_i$ and ending in vertex $v_j$ in $\Delta$.

For each $n \in \N$ take a finite full subquiver $\Delta_n$ of $\Delta$ (say, on the vertices $\{ v_{n,1},\hdots,v_{n,k} \}$) then \cite{AdriLieven} asserts that there is an \'etale map between a Zariski open subset of the moduli space $M_{\alpha}^{ss}(Q,\theta)$ of $\theta$-semi-stable $Q$-representations of dimension vector $\alpha = \alpha_{n,1} + \alpha_{n,2} + \hdots + \alpha_{n,k}$ around the $Q$-representation $V_{n,1} \oplus V_{n,2} \oplus \hdots \oplus V_{n,k}$ {\em and} the moduli space of {\em semi-simple} $\Delta_n$-representations of dimension vector $\pmb{1} = (1,1,\hdots,1)$ around the zero-representation. Moreover, in this \'etale correspondence, (isomorphism classes of) simple $\Delta_n$-representations correspond to (isomorphism classes) of $\theta$-stable representations.

By the results from \cite{LBProcesi} we have accomplished our objective, provided we can find for each $n$ a {\em subquiver} $\Sigma_n$ of $\Delta_n$ satisfying the following conditions
\begin{itemize}
\item{$\Sigma_n$ is strongly connected, meaning that any two vertices are connected via an oriented circuit in $\Sigma_n$, and}
\item{$1-\chi_{\Sigma_n}(\pmb{1},\pmb{1}) = k(n)-1$}
where $\chi_{\Sigma_n}$ is the Euler-form (as above) of the quiver $\Sigma_n$.
\end{itemize}
An example : consider the set $\{ V_0 = S_{11}, V_1=T_1(\lambda_1),V_2=T_2(\lambda_2),V_3=T_1(\lambda_3),V_4=T_2(\lambda_4),V_5=T_1(\lambda_5),\hdots \}$ with $\lambda_i \not= \lambda_j$ if $i \not= j$. Then, the quiver $\Delta$ has exactly one loop in each vertex $v_i$ (except in $v_0$)  and exactly one arrow $v_i \rTo v_j$ whenever $i \not=j~mod~2$. Let $\Delta_n$ be the full subquiver on the first $[ \tfrac{n}{2} ]$ vertices and $\Sigma_n$ the subquiver below (on vertices $\{ v_1,\hdots, v_{[ \tfrac{n}{2} ]}\}$ if $n$ is even and on $\{ v_0,v_1,\hdots,v_{[ \tfrac{n}{2} ]} \}$ if $n$ is odd). Then, the indicated representations give an $[ \tfrac{n}{2} ]$-parameter family of simple $\Sigma_n$ (and hence also $\Delta_n$)-representations
\[
\text{($n$ even) : }~\xymatrix{
\vtx{1} \ar@(ul,dl)|{\alpha_1} \ar@/^/[r]|{\alpha_2} & \vtx{1} \ar@/^/[l]|{1} \ar@/^/[r]|{\alpha_3} & \vtx{1} \ar@/^/[r]|{\alpha_4} \ar@/^/[l]|{1} & \ar@/^/[l]|{1} \hdots & \hdots \ar@/^/[r] & \vtx{1} \ar@/^/[r]|{\alpha_{[ \tfrac{n}{2} ]}} \ar@/^/[l] & \vtx{1} \ar@/^/[l]|{1} }
\]
\[
\text{($n$ odd) : }~\xymatrix{
\vtx{1} \ar@/^/[r]|{\alpha_1} & \vtx{1}  \ar@/^/[r]|{\alpha_2} \ar@/^/[l]|{1} & \vtx{1} \ar@/^/[l]|{1} \ar@/^/[r]|{\alpha_3} & \vtx{1} \ar@/^/[r]|{\alpha_4} \ar@/^/[l]|{1} & \ar@/^/[l]|{1} \hdots & \hdots \ar@/^/[r] & \vtx{1} \ar@/^/[r]|{\alpha_{[ \tfrac{n}{2} ]}} \ar@/^/[l] & \vtx{1} \ar@/^/[l]|{1} }
\]
Using the \'etale map these representations give an $[ \tfrac{n}{2} ]$-parameter family of $\theta$-stable $Q$-representations and hence of irreducible $n$-dimensional $\Gamma_0$-representations, and hence by Schur an $[ \tfrac{n}{2} ] + 1$-parameter family of isomorphism classes of irreducible $B_3$-representations.

\par \vskip 4mm

\noindent
{\bf (step 4) Reverse-engineer} the above general argument to fit your specific example.

 \end{document}